\documentclass[12pt]{article}
\usepackage{amsthm,amssymb,amsmath,graphicx,cite}
\usepackage{algorithmic, algorithm}
\usepackage[margin = 1in]{geometry}

\newcommand{\LHS}{\mathsf{LHS}}

\newcommand{\Oh}{{\mathcal O}}
\newcommand{\dmin}{\displaystyle\min}

\newcommand{\ip}[2]{\left\langle #1 , #2 \right\rangle}    
\newcommand{\R}{{\mathbb R}}

\newcommand{\prox}{{\text{\sf Prox}}}

\newcommand{\dom}{{\mathrm{dom}}}

\DeclareMathOperator*{\argmin}{arg\,min}

\DeclareMathOperator{\dist}{dist}

\DeclareMathOperator{\interior}{int}

\newtheorem{lemma}{Lemma}
\newtheorem{theorem}{Theorem}

\title{Convergence of first-order methods via the convex conjugate}
\author{
Javier Pe\~na\thanks{Tepper School of Business,
Carnegie Mellon University, USA, {\tt jfp@andrew.cmu.edu}}
}

\begin{document}

\maketitle

\begin{abstract}
This paper gives a unified and succinct approach to the $\Oh(1/\sqrt{k}), \Oh(1/k),$ and $\Oh(1/k^2)$ convergence rates of the subgradient, gradient, and accelerated gradient methods for unconstrained convex minimization.  In the three cases the proof of convergence follows from  a generic bound defined by the convex conjugate of the objective function.
\end{abstract}


\section{Introduction}
\label{sec.intro}

The subgradient, gradient, and accelerated gradient methods are icons in the class of first-order algorithms for convex optimization.  Under a suitable Lipschitz continuity assumption on the objective function and a judicious choice of step-sizes, the subgradient method yields a point whose objective value is within $\Oh(1/\sqrt{k})$ of the optimal value after $k$ iterations.  In a similar vein, under a suitable Lipschitz continuity assumption on the gradient of the objective function and a judicious choice of step-sizes, the gradient and accelerated gradient methods yield points whose objective values are within $\Oh(1/k)$ and $\Oh(1/k^2)$ of the optimal value respectively after $k$ iterations.  

Although the proofs of the $\Oh(1/\sqrt{k}), \Oh(1/k),$ and $\Oh(1/k^2)$ convergence rates for these three algorithms share some common ideas, they are traditionally treated separately.  
In particular, the known proofs of the $\Oh(1/k^2)$ convergence rate of the accelerated gradient method, first established by Nesterov in a landmark paper~\cite{Nest83}, are notoriously less intuitive than those of the $\Oh(1/\sqrt{k})$ and $\Oh(1/k)$ convergence rates of the subgradient and gradient methods.  Nesterov's accelerated gradient method has had a profound influence in optimization and has led to a vast range of developments.  See, e.g.,~\cite{BeckCG11,BeckT09,SraNW12,Nest04,Tsen08} and the many references therein.

Several recent articles~\cite{AlleO14,BubeLS15,FlamB15,LessRP16,ODonC15,SuBC14} have proposed novel approaches that add insight and explain how the accelerated gradient method and some variants achieve a faster convergence rate.  This paper makes a  contribution of similar spirit.  It provides a unified and succinct approach for deriving the convergence rates of the subgradient, gradient, and accelerated gradient algorithms.  The crux of the approach is a generic upper bound via the convex conjugate of the objective function. (See Lemma~\ref{the.lemma} in Section~\ref{sec.main}.)  The construction of the upper bound captures key common features and differences among the three algorithms. 

\medskip

The paper is self-contained and relies only  on the basic convex analysis background  recalled next.  (For further details see~\cite{BorwL00,HiriL93,Rock70}.)  Let $f:\R^n \rightarrow \R\cup\{\infty\}$ be a convex function.  Endow $\R^n$ with an inner product $\ip{\cdot}{\cdot}$ and let $\|\cdot\|$ denote the corresponding Euclidean norm.  Given a constant $G >0$, the function $f$ is $G$-Lipschitz if  for all $x,y\in \dom(f):= \{x\in \R^n: f(x) < \infty\}$
\[f(x) - f(y) \le G\|x-y\|.\]
Observe that if $f$ is convex and $G$-Lipschitz then for all $x\in \interior(\dom(f))$ and $g\in \partial f(x)$
\begin{equation}\label{eq.Lip.fun}
g \in \partial f(x) \Rightarrow \|g\| \le G.
\end{equation}
Suppose $f$ is differentiable on $\dom(f)$.  Given a constant $L>0$, the gradient $\nabla f$ is $L$-Lipschitz if  for all $x,y\in \dom(f)$
\[\|\nabla f(x) - \nabla f(y)\| \le L\|x-y\|.\]
Observe that if $f$ is differentiable and $\nabla f$ is $L$-Lipschitz then for all $x,y\in \dom(f)$
\[f(y) \le f(x) + \ip{\nabla f(x)}{y-x} + \frac{L}{2}\|y-x\|^2.\]   
In particular, if $x\in \dom(f)$ is such that $x-\frac{1}{L}\nabla f(x) \in \dom(f)$ then
\begin{equation}\label{eq.Lip}
f\left(x-\frac{1}{L}\nabla f(x)\right) \le f(x) - \frac{1}{2L}\|\nabla f(x)\|^2.
\end{equation}
Let $f^*:\R^n\rightarrow \R\cup \{\infty\}$ denote the convex conjugate of $f$, that is,
\[
f^*(z) = \sup_{x\in \R^n} \left\{\ip{z}{x} - f(x)\right\}.
\]
The construction of the conjugate readily yields the following property known as {\em Fenchel's inequality.} For all $z,x\in \R^n$
\[
f^*(z) + f(x) \ge \ip{z}{x}
\]
and equality holds if $z \in \partial f(x)$. 

\section{First-order methods for unconstrained convex optimization}
\label{sec.main}
Throughout the sequel assume $f:\R^n \rightarrow \R$ is a convex function and consider the problem 
\begin{equation}\label{eq.problem}
\dmin_{x\in \R^n}  f(x). 
\end{equation} 
Let $\bar f$ and $\bar X$ respectively denote the optimal value and set of optimal solutions to \eqref{eq.problem}.  

Algorithm~\ref{algo.subgrad} and Algorithm~\ref{algo.accel.grad}  describe respectively the subgradient method and accelerated gradient method for \eqref{eq.problem}.  The subgradient method becomes the gradient method when $f$ is differentiable.  Algorithm~\ref{algo.accel.grad} is a variant of Nesterov's original accelerated gradient method~\cite{Nest83}.  This version has been discussed in~\cite{BeckT09,Nest04,Tsen08}.  

\begin{algorithm}
\caption{Subgradient/gradient method}\label{algo.subgrad}
\begin{algorithmic}[1]
	\STATE {\bf input:} $x_0 \in \R^{n}$ and a convex  function $f:\R^n \rightarrow \R$
	\FOR{$k=0,1,2,\dots$}
	 	\STATE pick $g_k \in \partial f(x_k)$ and $t_k > 0$
		\STATE $x_{k+1}:=x_k-t_k g_k$
	\ENDFOR
\end{algorithmic}
\end{algorithm}

\begin{algorithm}
\caption{Accelerated gradient method}\label{algo.accel.grad}
\begin{algorithmic}[1]
	\STATE {\bf input:} $x_0 \in \R^{n}$  and a differentiable convex  function $f:\R^n \rightarrow \R$
	\STATE $y_0:=x_0,$ $\theta_0:=1$
	\FOR{$k=0,1,2,\dots$}
		\STATE pick $t_k > 0$
		\STATE $x_{k+1}:=y_k-t_k\nabla f(y_k)$
		\STATE let $\theta_{k+1} \in (0,1)$ be such that $\theta_{k+1}^2 = \theta_k^2(1-\theta_{k+1})$
		\STATE $y_{k+1}:=x_{k+1} + \frac{\theta_{k+1}(1-\theta_k)}{\theta_k} (x_{k+1}-x_k)$
	\ENDFOR
\end{algorithmic}
\end{algorithm}
Theorem~\ref{thm.subgrad}, Theorem~\ref{thm.grad}, and Theorem~\ref{main.thm} state well-known convergence properties of Algorithm~\ref{algo.subgrad} and Algorithm~\ref{algo.accel.grad}.

\begin{theorem}\label{thm.subgrad} Suppose $f$ is $G$-Lipschitz.  Then the sequence of iterates $x_k \in \R^n, \; k=1,2,\dots$ generated by Algorithm~\ref{algo.subgrad} satisfies
\begin{equation}\label{eq.thm.subgrad}
\frac{ \sum_{i=0}^{k} t_i f(x_i) - \frac{G^2}{2} \sum_{i=0}^{k} t_i^2}{\sum_{i=0}^{k} t_i}
\le   f(x) +\frac{\|x_0- x\|^2}{2\sum_{i=0}^{k} t_i}
\end{equation}
for all $ x \in \R^n.$ In particular, if $\bar X\ne \emptyset$ then
${\dmin_{i=0,1,\dots,k} f(x_i)} - \bar f \le \frac{\dist(x_0,\bar X)^2 + G^2 \sum_{i=0}^{k} t_i^2}{2\sum_{i=0}^{k} t_i},$
and 
${\dmin_{i=0,1,\dots,k} f(x_i)} - \bar f \le \frac{\dist(x_0,\bar X)^2 + G^2}{2\sqrt{k+1}}$
for $t_i = \frac{1}{\sqrt{k+1}}, \; i=0,1,\dots,k.$
\end{theorem}

\begin{theorem}\label{thm.grad} Suppose $\nabla f$ is $L$-Lipschitz and $t_k = \frac{1}{L}, \; k=0,1,\dots$. Then the sequence of iterates $x_k \in \R^n, \; k=1,2,\dots$ generated by Algorithm~\ref{algo.subgrad} satisfies
\begin{equation}\label{eq.thm.grad}
\frac{f(x_1) + \cdots + f(x_k)}{k} \le f(x) + \frac{L \|x_0-x\|^2}{2k}
\end{equation}
for all $x \in \R^n.$  In particular,  if $\bar X\ne \emptyset$ then
$f(x_k) - \bar f \le \frac{L\dist(x_0,\bar X)^2}{2k}.$
\end{theorem}

\begin{theorem}\label{main.thm} Suppose $f$ is differentiable, $\nabla f$ is $L$-Lipschitz, and $t_k = \frac{1}{L}$  for $k=0,1,\dots$. Then the sequence of iterates $x_k \in \R^n, \; k=1,2,\dots$ generated by Algorithm~\ref{algo.accel.grad} satisfies
\begin{equation}\label{eq.main.thm}
f(x_k) \le f(x) + \frac{L \theta_{k-1}^2 \|x_0-x\|^2}{2}
\end{equation}
for all $x \in \R^n.$ In particular, if $\bar X\ne \emptyset$ then
$f(x_k) - \bar f \le \frac{2L \dist(x_0,\bar X)^2}{(k+1)^2}.$
\end{theorem}
The central contribution of this paper is a unified  approach to  the proofs of Theorem~\ref{thm.subgrad}, Theorem~\ref{thm.grad}, and Theorem~\ref{main.thm}.  The crux of the approach is the following lemma.

\begin{lemma}\label{the.lemma}  There exists a sequence $z_k\in \R^n,\; k=1,2,\dots$ such that for $k=1,\dots$ and $\mu_k = \frac{1}{\sum_{i=0}^k t_i}$ the left-hand side $\LHS_k$ of \eqref{eq.thm.subgrad} in Theorem~\ref{thm.subgrad}  satisfies
\begin{equation}\label{eq.lemma}
\LHS_k \le  -f^*(z_k) + \ip{z_k}{x_0} - \frac{\|z_k\|^2}{2\mu_k} = -f^*(z_k) + \min_{u\in\R^n}\left\{\ip{z_k}{u} + \frac{\mu_k}{2}\|u-x_0\|^2  \right\}.
\end{equation}
There also exist sequences $z_k\in \R^n,\; k=1,2,\dots$ such that~\eqref{eq.lemma} holds for $\mu_k = \frac{L}{k}$ and the
left-hand side $\LHS_k$ of~\eqref{eq.thm.grad} in Theorem~\ref{thm.grad}, as well as for~$\mu_k = L\theta_{k-1}^2$ and the
left-hand side $\LHS_k$ of~\eqref{eq.main.thm} in Theorem~\ref{main.thm}.
\end{lemma}

\medskip

Lemma~\ref{the.lemma} captures some key common features and differences among the subgradient, gradient, and accelerated gradient algorithms.  The right-hand side in \eqref{eq.lemma} has the same form in all cases and has the same kind of dependence on the initial point $x_0$. Furthermore, as Section~\ref{sec.proof} below details, the construction of the sequences $z_k,\mu_k,  \; k = 1,2\dots$ follows the same template for the three algorithms.  However, some details of the construction  for these sequences need to be carefully tailored to each of the three algorithms.

\medskip

\noindent
{\bf Proof of Theorem~\ref{thm.subgrad}, Theorem~\ref{thm.grad}, and Theorem~\ref{main.thm}.}
Lemma~\ref{the.lemma} and Fenchel's inequality imply that for some $z_k\in \R^n,\;k=1,2,\dots$ and 
all $ x \in \R^n$  the left-hand-sides $\LHS_k$ of~\eqref{eq.thm.subgrad}, \eqref{eq.thm.grad}, and \eqref{eq.main.thm} satisfy
\begin{align*}
\LHS_k&\le  -f^*(z_k)  + \min_{u\in\R^n}\left\{\ip{z_k}{u} + \frac{\mu_k}{2}\|u-x_0\|^2  \right\}  \\
 & \le -f^*(z_k) + \ip{z_k}{x} + \frac{\mu_k \cdot \|x - x_0\|^2}{2}\\
 & \le   f(x) + \frac{\mu_k \cdot \|x - x_0\|^2}{2}.
\end{align*}
To finish, recall that $\mu_k = \frac{1}{\sum_{i=0}^k t_i}$ for \eqref{eq.thm.subgrad}, $\mu_k = \frac{L}{k}$ for~\eqref{eq.thm.grad}, and $\mu_k = L\theta_{k-1}^2$ for~\eqref{eq.main.thm}.
For the second part of Theorem~\ref{thm.grad} observe that $f(x_k) \le \frac{f(x_1) + \cdots + f(x_k)}{k}$ because~\eqref{eq.Lip} implies that $f(x_{i+1}) \le f(x_i) - \frac{1}{2L}\|\nabla f(x_i)\|^2 \le f(x_i),\;i=0,1,\dots.$
For the second part of Theorem~\ref{main.thm} observe that a straightforward induction shows that the conditions 
$\theta_{k+1} \in (0,1), \theta_{k+1}^2 = \theta_k^2(1-\theta_{k+1}),$ and $\theta_0 = 1$  imply $\theta_{k-1} \le \frac{2}{k+1}$.
\qed

\section{Proof of Lemma~\ref{the.lemma}}
\label{sec.proof}
Construct the sequences $\mu_k\in \R, \; z_k\in \R^n, \; k = 1,2\dots$ as follows.  First, choose sequences $\theta_k \in (0,1), y_k \in \R^n, g_k \in \partial f(y_k),\; k = 1,2,\dots,$ and two initial values  $\mu_0 \in \R_+, z_0\in \R^n$ or $\mu_1 \in \R_+, z_1\in \R^n$.  Second, let $\mu_k\in \R, \; z_k\in \R^n, \; k = 1,2\dots$ be defined by the rules
\[
\begin{array}{rcl}
z_{k+1} &=& (1-\theta_k) z_k + \theta_k g_k \\
\mu_{k+1}&=& (1-\theta_k) \mu_k.
\end{array}
\]
This construction readily implies 
\[
\begin{array}{rcl}
\ip{z_{k+1}}{x_0} - \frac{\|z_{k+1}\|^2}{2\mu_{k+1}} &=& (1-\theta_k) \left(\ip{z_k}{x_0} - \frac{\|z_{k}\|^2}{2\mu_{k}}\right) \\
&& + \theta_k \left(\ip{g_k}{x_0-\frac{z_k}{\mu_k}} -\frac{\theta_k}{2(1-\theta_k) \mu_k} \|g_k\|^2\right),
\end{array}
\]
and, by the convexity of $f^*$ and $g_k \in \partial f(y_k)$,
\[
\begin{array}{rcl}
-f^*(z_{k+1}) &\ge& -(1-\theta_k)f^*(z_k) - \theta_k f^*(g_k) \\&=& -(1-\theta_k)f^*(z_k) - \theta_k \left(\ip{g_k}{y_k} + f(y_k)\right).
\end{array}
\]
Thus
\begin{multline}\label{eq.induction}-f^*(z_{k+1}) + \ip{z_{k+1}}{x_0} - \frac{\|z_{k+1}\|^2}{2\mu_{k+1}} \ge (1-\theta_k)\left(-f^*(z_k) +\ip{z_k}{x_0} - \frac{\|z_{k}\|^2}{2\mu_{k}}\right) \\ + \theta_k \left(\ip{g_k}{x_0-y_k-\frac{z_k}{\mu_k}} + f(y_k)-\frac{\theta_k}{2(1-\theta_k) \mu_k} \|g_k\|^2\right).
\end{multline}
To prove~\eqref{eq.lemma}, proceed  by induction.  By~\eqref{eq.induction} to show the inductive step $k$ to $k+1$ it suffices to show
\begin{equation}\label{eq.to.show}
\LHS_{k+1} - (1-\theta_k)\LHS_k \le \theta_k \left(\ip{g_k}{x_0-y_k-\frac{z_k}{\mu_k}} + f(y_k)-\frac{\theta_k}{2(1-\theta_k) \mu_k} \|g_k\|^2\right).
\end{equation}
Next show \eqref{eq.to.show} in each of the three cases.  

First, for~\eqref{eq.thm.subgrad} take $\theta_k = \frac{t_{k+1}}{\sum_{i=0}^{k+1} t_i}, y_k = x_{k+1},$ and initial values $\mu_0 = \frac{1}{t_0}, z_0 = \frac{t_0 g_0}{t_0} = g_0$.  Then $\mu_k = \frac{1}{\sum_{i=0}^k t_i}, \; \frac{\theta_k}{(1-\theta_k)\mu_k} = t_{k+1},$ and $x_0-y_k - \frac{z_k}{\mu_k} = 0$.  Therefore
\[
\begin{array}{rcl}
\LHS_{k+1} - (1-\theta_k)\LHS_k &=& \frac{t_{k+1} f(x_{k+1}) - \frac{G^2}{2} t_{k+1}^2}{\sum_{i=0}^{k+1} t_i} \\
&=& 
\theta_k\left(f(x_{k+1}) - \frac{\theta_k}{2(1-\theta_k) \mu_k} G^2\right) 
\\&\le& \theta_k \left(\ip{g_k}{x_0-y_k-\frac{z_k}{\mu_k}} + f(y_k)-\frac{\theta_k}{2(1-\theta_k) \mu_k} \|g_k\|^2\right).
\end{array}
\]
The inequality in the last step follows from~\eqref{eq.Lip.fun}.

Second, for~\eqref{eq.thm.grad} take $\theta_k = \frac{1}{k+1}, y_k=x_{k},$ and initial values $\mu_1 = L, z_1 = \nabla f(x_0)$.  Then $\mu_k = \frac{L}{k},\; \frac{\theta_k}{(1-\theta_k)\mu_k} = L,$ and $x_0-y_k - \frac{z_k}{\mu_k} = 0$.  Therefore
\[
\begin{array}{rcl}
\LHS_{k+1} - (1-\theta_k)\LHS_k &=& \frac{f(x_{k+1})}{k+1} \\
&\le& \theta_k(f(x_{k})-\frac{1}{2L}\|\nabla f(x_k)\|^2)
\\
&=& 
\theta_k\left(f(y_{k}) - \frac{\theta_k}{2(1-\theta_k) \mu_k} \|g_k\|^2\right) 
\\&=& \theta_k \left(\ip{g_k}{x_0-y_k-\frac{z_k}{\mu_k}} + f(y_k)-\frac{\theta_k}{2(1-\theta_k) \mu_k} \|g_k\|^2\right).
\end{array}
\]
The inequality in the second step follows from $x_{k+1} = x_k - \frac{1}{L}\nabla f(x_k)$ and~\eqref{eq.Lip}.

Third, for~\eqref{eq.main.thm} take $\theta_k, y_k$ as in Algorithm~\ref{algo.accel.grad} and initial values $\mu_1 = L, z_1=\nabla f(x_0)$.  A separate  induction argument shows that $\mu_{k} = L \theta_{k-1}^2,\; \frac{\theta_k^2}{(1-\theta_k)\mu_k} = \frac{1}{L},$ and 
\[
\begin{array}{rcl}
y_{k}&=&(1-\theta_k)x_k + \theta_k(x_0-\frac{z_k}{\mu_k})\\
x_{k+1}&=&(1-\theta_k)x_k + \theta_k(x_0-\frac{z_{k+1}}{\mu_{k+1}})
\end{array}
\]
for $k=1,2,\dots.$ In particular,
\begin{equation}\label{eq.kstep}
(1-\theta_k)(y_k-x_k) = \theta_k\left(x_0 - y_k - \frac{z_k}{\mu_k}\right).
\end{equation}
Therefore
\[
\begin{array}{rcl}
\LHS_{k+1} - (1-\theta_k)\LHS_k &=& f(x_{k+1}) - (1-\theta_k)f(x_k) \\
&\le& 
f(y_k) -\frac{1}{2L}\|\nabla f(y_k)\|^2 - (1-\theta_k)(f(y_k) + \ip{\nabla f(y_k)}{x_k-y_k})
\\
&=& 
(1-\theta_k)\ip{g_k}{y_{k}-x_{k}} + \theta_k f(y_k)  - \frac{1}{2L}\|g_k\|^2 
\\
&= & \theta_k \left(\ip{g_k}{x_0-y_k-\frac{z_k}{\mu_k}} + f(y_k)-\frac{\theta_k}{2(1-\theta_k) \mu_k} \|g_k\|^2\right).
\end{array}
\] 
The inequality in the second step follows from $x_{k+1} = y_k - \frac{1}{L}\nabla f(y_k)$ and~\eqref{eq.Lip}, and from the convexity of $f$. The fourth step follows from 
\eqref{eq.kstep}. 

To complete the proof of~\eqref{eq.lemma} by induction it only remains to verify that \eqref{eq.lemma} holds for $k=0$ or $k=1$ in each of the three cases.  For~\eqref{eq.thm.subgrad} observe that 
$
f(x_0) =  \ip{z_0}{x_0} - f^*(z_0)
$ because $z_0 = g_0\in \partial f(x_0)$.  From~\eqref{eq.Lip.fun} and $\mu_0 = \frac{1}{t_0}$ it follows that
\[
\LHS_0 = \frac{t_0f(x_0) - \frac{G^2}{2}t_0^2}{t_0} = 
f(x_0) - \frac{t_0}{2}G^2 \le -f^*(z_0) + \ip{z_0}{x_0} -\frac{\|z_0\|^2}{2\mu_0}.
\]
For both~\eqref{eq.thm.grad} and~\eqref{eq.main.thm} observe that 
$
f(x_0) =  \ip{z_1}{x_0} - f^*(z_1)
$
because $z_1 = \nabla f(x_0)$. From~\eqref{eq.Lip} and $\mu_1 = L$, it follows that
\[
\LHS_1 = f(x_1) = f\left(x_0 - \frac{1}{L}\nabla f(x_0)\right) \le f(x_0) - \frac{1}{2L} \|\nabla f(x_0)\|^2 = -f^*(z_1) + \ip{z_1}{x_0} -\frac{\|z_1\|^2}{2\mu_1}.
\]

\section{Potential extensions}
This section sketches some potential  extensions that will a topic for future work. 

\subsection{Proximal iterations}
There are various first-order methods defined via  proximal iterations~\cite{BeckCG11,BeckT09,ChamD15,Gule91,Tsen08}.  Suppose $f=\phi + \psi$, where $\phi, \psi:\R^n \rightarrow \R\cup\{\infty\}$ are convex functions such that the proximal map
\[
\prox_{t}(x):= \argmin_{y\in \R^n}\left\{\psi(y) + \frac{1}{2t}\|x-y\|^2\right\}
\]
is computable.  If $\phi$ is differentiable, then Algorithm~\ref{algo.accel.grad}  extends (see, e.g.,~\cite{BeckT09}) by replacing step 5 with 
\[
5': \; x_{k+1} = \prox_{t_k}(y_k - t_k \nabla \phi(y_k)).
\]
Algorithm~\ref{algo.subgrad} also extends in a similar fashion.  

A suitable extended version of Lemma~\ref{the.lemma} would readily yield a unified proof of the corresponding extended versions of Theorem~\ref{thm.subgrad}, Theorem~\ref{thm.grad}, and Theorem~\ref{main.thm}.  The author conjectures that this is indeed the case if the right hand side in~\eqref{eq.lemma} is replaced with the following expression
\[
-\phi^*(z_k) + \min_{u\in\R^n}\left\{\psi(u)+\ip{z_k}{u} + \frac{\mu_k}{2}\|u-x_0\|^2  \right\}.
\] 

\subsection{Stronger convergence results}

The convex conjugate approach developed in this paper may also yield alternative proofs of other stronger  convergence properties of first-order methods.  In particular, the $\Oh(1/k)$ and $\Oh(1/k^2)$ convergence rates of the gradient and the accelerated gradient methods the can be strengthened to $o(1/k)$ and $o(1/k^2)$ respectively as shown in~\cite{AttoP16,Gule91}.  It is also known that the sequence of iterates generated by the gradient and accelerated gradient methods converge weakly to a minimizer as discussed in~\cite{AttoCPR16,ChamD15}.  The convex conjugate approach introduced in this paper  may lead to succinct and unified derivations of these and possibly other results.  

\section*{Acknowledgements}

This research has been  funded by NSF grant CMMI-1534850.

\end{document}